# Probable Counterexamples of the *ABC* Conjecture
N. A. Carella, June 2003

*Abstract:* This paper describes a method used to construct infinitely many probable counterexamples of the *abc* conjecture over the rational integers.

**Introduction**
The *abc* conjecture, due to Masser and Oesterle, circa 1985, is a statement in arithmetical geometry about the integer solutions of the equation $x + y = z$ in relatively prime triples *a*, *b*, *c*. The precise statement is as follows.

***Conjecture* 1.** (The *abc* conjecture)   Given a real number $\varepsilon > 0$, there exists a constant $c_\varepsilon$ such that if $a + b = c$, and $\gcd(a, b, c) = 1$, then

$$\max\{|a|, |b|, |c|\} \leq c_\varepsilon \prod_{p | abc} p^{1+\varepsilon}. \qquad (1)$$

There are a few versions of it over the integers, and generalizations to numbers fields, and functions fields, et cetera. The corresponding claim in ring of polynomial functions, known as Mason's theorem, is easy to prove, see [7], [9].

Every integer $n \in \mathbb{N}$ has exactly $[n/2]$ partitions as a sum of two integers, and exactly $\varphi(n)/2$ partitions as a sum of two relatively prime integers, some details on the maximal number of partitions that satisfy (1) appear in [16].

The radical of an integer *n* is defined by the expression $rad(n) = \prod_{p|n} p$, and the radical index is defined by ratio $\nu(n) = \log(n) / \log(rad(n))$. A highly composite (smooth) integer has a large radical index. The conjecture states that the sum (or difference) of two relatively prime integers can not have a large radical index; it is less than some absolute constant. This ratio is currently used to test for bad cases of the conjecture. A few authors are tabulating high values of the compositeness ratio. The record is 1.6299 for the triple $2 + 109 \cdot 3^{10} = 23^5$, see [12]. In [1] there is a discussion about a method (using continued fractions) for generating triples with high ratios.

On the last two decades this conjecture has been the object of extensive research efforts. Some of the techniques developed as result of these efforts are relatively easy to work with. Accordingly the analysis of some problems (some still unsolved) is conceptually easy to handle with the *abc* conjecture. For these reasons the applications of the *abc* conjecture are growing. Presently the applications range from the solvability of algebraic equations, powers in integer recurring sequences, parameters of elliptic curves and modular symbols, to the cardinality of the set of Wieferich primes, and numerical analysis, see [2], [9], [11], [17], [6], [12], etc; the latter has a list of about 26 different applications. In light of the result herein, some of these abc conditional results possibly will have to be reassessed.

**Preliminary**

The ensuing results are under the context of an algebraic curve $C : y^2 = x^3 + ax + b$ and its group of rational points $C(\mathbb{Q}) \cong T_{tor} \times \mathbb{Z}^r$ of rank $r$ and torsion subgroup $T_{tor}$. An arbitrary point $P \in C(\mathbb{Q}) = \{ n_r P_r + \cdots + n_1 P_1 + T : n_i \in \mathbb{Z} \}$ is a linear combination of the generators and the torsion points. A multiple of a point is written as $nP = (x(nP), y(nP))$.

**Lemma 2.** Let $x(nP) = x_{nP}/z_{nP}^2$, then $z_{nP}^2 = x_{nP} e^{-\alpha(n)}$, where $\alpha(n) = o(n^2)$.

Proof: It is known that the asymptotic limit of the ratio

$$\frac{\log |x_{nP}|}{\log |z_{nP}^2|} \tag{2}$$

of the numbers of digits in the numerator and denominator of $x(nP)$ tends to unity. And the height of the point $nP = (x(nP), y(nP))$ satisfies the relation

$$h(nP) = \max\{ \log |x_{nP}|, \log |z_{np}^2| \} = 2n^2 h(P) + O(1), \tag{3}$$

where the constant depends only on the curve, see [4], and [10]. The claim follows from these data. ∎

**Definition 3.** A sequence of points $\{ P_n = (x_n/z_n^2, y_n/z_n^3) \} \subset C(\mathbb{Q})$ is said to have exponential growth if $x_n/z_n^2 = x_n^{\gamma(n)} \geq x_n^\delta$ for some $\delta > 0$ and all $n \geq 1$. Otherwise $\gamma(n) \to 0$ as $n \to \infty$, and it has subexponential growth.

Naturally different sequences of points can have different rates of growth. For example, if $r \geq 1$, the sequence of points $\{ nP : n \in \mathbb{N} \}$ has subexponential growth, this follows from Lemma 2. But for $r > 1$ a sequence of points $\{ P_n = nP + mQ : n \in \mathbb{N} \}$ can have exponential growth. Sequences of points on the curve $C : y^2 = x^3 + ax + b$ with exponential growth contradict the abc conjecture.

**Description of the Probable Counterexamples**

Let $\varepsilon = 0$, and consider the infinite sequence of triples

$$a_n = p^{q^{n-1}(q-1)} - 1, \quad b_n = 1, \quad c_n = a_n + b_n, \tag{4}$$

where $q$ is prime and $\gcd(p, q) = 1$, see [9]. Since $q^n \mid p^{q^{n-1}(q-1)} - 1$ for all $n > 0$, (this follows from Fermat's little theorem), one concludes that for any fixed constant $c_0$, inequality (1) fails infinitely often as $n$ or $q \to \infty$, $p$ fixed. This confirms the requirement that $\varepsilon \neq 0$. The authors in [18] have proven the existence of infinitely many integers such that

$$\max\{|a|, |b|, |c|\} > N e^{(4-\delta)\sqrt{\log N / \log \log N}}, \tag{5}$$

where $N = \operatorname{rad}(abc)$, $\delta > 0$. This confirms the requirement that $\varepsilon > (4-\delta)/\sqrt{\log N \log \log N}$ as $N \to \infty$, and that $\varepsilon$ is not a rapidly decreasing function of $N$. Their work was later slightly improved in [19]. On the other direction there is the estimate

$$a + b = c < e^{c_1 N^{1/3} (\log N)^3}, \tag{6}$$

where $c_1$ is an absolute constant, see [19].

A few authors claim that there is no hope of proving or disproving the *abc* conjecture, and other claim that there is hope of disproving it.

*Claim:* For any real numbers $\varepsilon > 0$, and a constant $c_\varepsilon$, it is plausible that there exists infinitely many triples of integers $c = a + b$, $\gcd(a, b, c) = 1$, such that

$$\max\{|a|, |b|, |c|\} \geq c_\varepsilon \prod_{p \mid abc} p^{1+\varepsilon} \tag{7}$$

could hold.

*Reason:* Consider an algebraic curve $C : y^2 = x^3 + d$ of $\operatorname{rank}(C(\mathbb{Q})) = r > 1$, $d$ fixed, and let

$$P = \left(\frac{x_P}{z_P^2}, \frac{y_P}{z_P^3}\right), Q = \left(\frac{x_Q}{z_Q^2}, \frac{y_Q}{z_Q^3}\right) \in C(\mathbb{Q}) = \{n_r P_r + \cdots + n_1 P_1 : n_i \in \mathbb{Z}\}, \tag{8}$$

where $x_P, y_P, z_P \in \mathbb{Z}$. The $x$-coordinate of the sum/difference of a pair of points is given by

$$x(P \pm Q) = \frac{(y_P z_Q^3 \mp y_Q z_P^3)^2 - (x_P z_Q^2 + x_Q z_P^2)(x_P z_Q^2 - x_Q z_P^2)^2}{(x_P z_Q^2 - x_Q z_P^2)^2 z_P^2 z_Q^2}, \tag{9}$$

and the $y$-coordinate of the sum/difference of a pair of points is given by

$$y(P \pm Q) =$$
$$\frac{(y_P z_Q^3 \mp y_Q z_P^3)(2x_Q z_P^2 + x_P z_Q^2)(x_P z_Q^2 - x_Q z_P^2)^2 - (y_P z_Q^3 - y_Q z_P^3) - y_Q z_P^3 (x_P z_Q^2 - x_Q z_P^2)^3}{(x_P z_Q^2 - x_Q z_P^2)^3 z_P^3 z_Q^3}$$

For each point $P \pm Q = (X/Z^2, Y/Z^3) \in C(\mathbb{Q})$, the triple $(X, Y, Z) \in \mathbb{Z}$ is an integer solution of the equation

$$X^3 = Y^2 + dZ^6. \tag{10}$$

Now estimating the radical rad($dXYZ$) and replacing them in the abc inequality returns

$$X^3 \leq c_\varepsilon \prod_{p|abc} p^{1+\varepsilon}$$
$$\leq c_\varepsilon (dXYZ)^{1+\varepsilon}. \tag{11}$$

Plugging in the dominant variable terms in the numerators of $x(P \pm Q)$ and $y(P \pm Q)$ yield

$$(x_P^3)^3 \leq c_\varepsilon \left( x_P^{15/2} (x_P z_Q^2 - x_Q z_P^2) z_P \right)^{1+\varepsilon}$$
$$\leq c_\varepsilon \left( x_P^8 (x_P z_Q^2 - x_Q z_P^2) \right)^{1+\varepsilon} \tag{12}$$

Since the set $C(\mathbb{Q})$ is dense in $C(\mathbb{R})$, see [3] and [14], each neighborhood of a fixed point $Q$ contains infinitely many rational points $P$. Therefore, as the point $P$ varies in a neighborhood of the fixed point $Q$, the term

$$Z = (x_P z_Q^2 - x_Q z_P^2) z_P z_Q \tag{13}$$

in the denominators of $x(P \pm Q)$ and $y(P \pm Q)$ takes small integer values infinitely often. Specifically the $x$-coordinate of a sum/difference $x(P \pm Q) = X/Z^2 \geq X^\delta$ has exponential with $\delta > 0$ a small constant. Consequently, as the height of the point $P = (X/Z^2, Y/Z^3)$ increases, the left side of (12) increases faster then the right side and the abc inequality is contradicted infinitely often. In fact a subsequence of the sequence of points $\{ P_{n,m} = nP - mQ) : m, n \in \mathbb{N} \}$ contradict (1) infinitely often since (13) vanishes whenever $m = n$.

This indicates that the asymptotic of the radical index $\log(X^3) / \log(rad(dXYZ)) > 1$.

*Remark:* Observe that the function $f : C(\mathbb{Q}) \to \mathbb{Q}$, defined by $f(P) = x(P - Q)$ has a pole at $P = Q$. And the sequence of points in a neighborhood of an $n$-torsion point $Q$ probably have special properties different from the properties of a sequence of points in a neighborhood of a nontorsion point $Q$.

**Statistical *ABC* Theory**
It is probably correct that almost all integers divisible by high powers of primes are not the sums (or differences) of two integers each divisible by high powers of primes, (the conjecture claims that all are not). This is exemplified by the tables of factored numbers, numerical experiments, integers recurring sequences, and the large collection of known Diophantine equations that have no integer solutions or have finitely many integer solutions. For example,

(1) $x^n + y^n = z^n$, $n \geq 3$, has no integers solutions,
(2) $ax^k + by^m = cz^n$, with $1/k + 1/m + 1/n < 1$ has finitely many integers solutions,
(3) $ax^k + by^m = 2^n$, has finitely many integers solutions.

More details on these topics appear in [2], [5], [6], [8], [10], [13], [15], [20], etc.

A problem of interest would be to develop a *statistical abc theory* comparable to the statistical theory of the prime divisors counting function $\omega(n)$ = #{ prime $p : p \mid n$ } developed by Ramanujan, Hardy, Erdos, Kac, and some resent authors. In this theory the values of the function $\omega(n)$ are normally distributed with mean loglog $n$, and standard deviation $(\text{loglog } n)^{1/2}$. Furthermore, as $x \to \infty$, the set of integers $n \leq x$ such that

$$\mid \omega(n) - \text{loglog } n \mid > (\text{loglog } n)^{1/2+\varepsilon} \qquad (14)$$

is set of density zero.

The statistical abc theory should answers questions of this type:

(1) Are the integers $c = a + b$, $\gcd(a, b, c) = 1$, that satisfy inequality (1) normally distributed as $c \to \infty$.
(2) Is the set of integers $c = a + b$, $\gcd(a, b, c) = 1$, that satisfy inequality (7) a set of density zero?

Undoubtedly the mean and standard deviation of such distribution function are functions of $1 + \varepsilon$, and as $1 + \varepsilon$ increases there are fewer integers that satisfy inequality (7). Perhaps there are none if $1 + \varepsilon \geq c_2$ = absolute constant.